\documentclass[11pt]{article}
\usepackage{amsmath}
\usepackage{amssymb}
\usepackage{amsthm}
\newcommand\beq{\begin{equation}}
\newcommand\eeq{\end{equation}}
\newcommand\bea{\begin{eqnarray}}
\newcommand\eea{\end{eqnarray}}
\newcommand\beast{\begin{eqnarray*}}
\newcommand\eeast{\end{eqnarray*}}
\newcommand\nonu{\nonumber}
\newcommand\dstyle\displaystyle
\newcommand\sa{\smallskipamount}
\newcommand\ma{\medskipamount}
\newcommand\ba{\bigskipamount}
\newcommand\sLP{\\[\sa]}
\newcommand\mLP{\\[\ma]}
\newcommand\bLP{\\[\ba]}
\newcommand\mPP{\\[\ma]\indent}

\newcommand\CC{\mathbb{C}}

\newcommand\FF{\mathbb{F}}

\newcommand\PP{\mathbb{P}}
\newcommand\QQ{\mathbb{Q}}
\newcommand\RR{\mathbb{R}}
\newcommand\ZZ{\mathbb{Z}}

\newcommand\FSN{{\cal N}}
\newcommand\FSO{{\cal O}}

\newcommand\FSU{{\cal U}}

\newcommand\al\alpha
\newcommand\be\beta
\newcommand\ga\gamma
\newcommand\de\delta
\newcommand\ep\varepsilon
\newcommand\ze\zeta
\newcommand\tha\theta
\newcommand\ka\kappa
\newcommand\la\lambda
\newcommand\si\sigma
\newcommand\om\omega

\newcommand\Om{\Omega}
\newcommand\half{\frac12}

\newcommand\thalf{\tfrac12}
\newcommand\iy\infty

\newcommand\union{\cup}

\newcommand\Zpos{\ZZ_{>0}}
\newcommand\Znonneg{\ZZ_{\ge0}}

\newcommand\const{{\rm const.}\,}
\newcommand{\qhyp}[5]{\,\mbox{}_{#1}\phi_{#2}\left(
  \genfrac{}{}{0pt}{}{#3}{#4};#5\right)}
\newcommand\LHS{left-hand side}
\newcommand\RHS{right-hand side}

\numberwithin{equation}{section}
\newtheorem{theorem}{Theorem}[section]

\newtheorem{corollary}[theorem]{Corollary}
\newtheorem{Definition}[theorem]{Definition}

\newtheorem{Remark}[theorem]{Remark}
\newenvironment{remark}{\begin{Remark}\rm}{\end{Remark}}
\newtheorem{Example}[theorem]{Example}

\newcommand{\bisub}[2]{\genfrac{}{}{0pt}{1}{#1}{#2}}
\begin{document}
\title{$LU$ factorizations, $q=0$ limits, and
$p$-adic interpretations of some $q$-hypergeometric
orthogonal polynomials}
\author{Tom H. Koornwinder and Uri Onn
\footnote{work done at KdV Institute, Amsterdam and supported by NWO,
project number 613.006.573}}
\date{Dedicated to Dick Askey on the occasion of his seventieth birthday}
\maketitle
\begin{abstract}
\noindent
For little $q$-Jacobi polynomials and $q$-Hahn polynomials
we give particular $q$-hypergeometric
series representations in which the termwise $q=0$ limit can be
taken. When rewritten in matrix form, these series representations
can be viewed as $LU$ factorizations. We develop a general theory of
$LU$ factorizations related to complete systems of orthogonal
polynomials with discrete orthogonality relations which admit a
dual system of orthogonal polynomials. For the $q=0$ orthogonal
limit functions we discuss interpretations on $p$-adic spaces. In
the little 0-Jacobi case we also discuss product formulas.
\end{abstract}
%
%until sec:08
%
\section{Introduction}
\label{sec:03}
This paper is concerned with limits for $q\downarrow0$ of some
$q$-hypergeometric orthogonal polynomials, in particular little
$q$-Jacobi polynomials and $q$-Hahn polynomials.
Limits of $q$-hypergeometric polynomials as $q\uparrow1$ are well-known,
see \cite[Chapter 5]{1}. Many ($q$-)hypergeometric orthogonal polynomials
have interpretations in connection with (quantum) group representations,
for instance as spherical or intertwining functions,
matrix elements, Clebsch-Gordan coefficients and
Racah coefficients, see for instance \cite{9} and \cite{10}.
Often, the $q\uparrow1$ limit of the polynomials corresponds to the
$q\uparrow1$ limit from the quantum group to the classical group.

Limits of $q$-hypergeometric polynomials for $q\downarrow0$
have been considered
for $q$-ultraspherical polynomials (see \cite[\S5]{14})
and for more general Askey-Wilson polynomials
(see \cite[pp.~26--28]{11} and references given there).
The limit functions have interpretations as spherical functions on
homogeneous trees (see references in \cite[p.~28]{11}) and on
infinite distance-transitive graphs (see \cite{12}).
Note that homogeneous trees are locally compact but noncompact
homogeneous spaces of the group $GL(2, \QQ_p)$ ($\QQ_p$ the field of
$p$-adic numbers).
No geometric explanation of this $q\downarrow0$ limit is known,
see also the discussion in \cite{13}.
Macdonald considered the $q\downarrow0$ limit of Macdonald polynomials,
both for root system $A_n$ (yielding Hall-Littlewood polynomials,
see \cite[Ch.~III]{16}) and for general root systems
(see \cite[\S10]{17}). The limit functions have interpretations as spherical
functions on a $p$-adic Lie group (see \cite{15}), in particular
in the $A_{d-1}$ case on $GL(d, \QQ_p)$ (see \cite[Ch.~V]{16}).

Haran \cite{3} considered limits for $q\downarrow0$ of $q$-Hahn
polynomials, little $q$-Jacobi polynomials and little $q$-Laguerre
polynomials. He gave interpretations of these limit functions as
spherical or intertwining functions on compact $p$-adic groups; in
particular in the {\em little $0$-Jacobi} case on the group
$GL(d,\ZZ_p)$ ($\ZZ_p$ the ring of $p$-adic integers).
In \cite{18} Haran's ideas were extended to the higher rank case and the
importance of the cellular basis, defined in \cite{4}, was emphasized.
All these interpretations are valid with respect to the larger family of
{\em $p$-adic fields}, or more generally {\em non-Archimedean fields}.
Altogether, $q$-special functions clearly play an important role as
interpolants between representations of groups over all local fields
($\RR$ and $\CC$ for $q=1$, and non-Archimedean local fields for
$q=0$).

The present paper gives new proofs of Haran's \cite{3}
$q\downarrow0$ limit results for little $q$-Jacobi polynomials and
$q$-Hahn polynomials (sections \ref{sec:01} and \ref{sec:08})
by starting with a series representation for these
polynomials which can be viewed as writing the square (possibly infinite)
matrix corresponding to the orthogonal polynomials as a product of
a lower triangular matrix and an upper triangular matrix
(an {\em $LU$ factorization}).
The matrix elements of the upper triangular matrix (simply $(q^x;q^{-1})_k$)
are the $q$-analogues of the cellular basis in the rank one case of \cite{4}.
This observation was decisive for the second author in order to find the
higher rank analogue
of this $q\downarrow0$ limit, see \cite{18}.
For the limiting little $0$-Jacobi functions we give a product formula.
In section \ref{sec:06} we discuss the $LU$ factorization more generally
for orthogonal polynomials with discrete orthogonality relations
which form a complete orthogonal system and for which
the dual orthogonal system also
consists of orthogonal polynomials.
We consider an upper-lower factorization as well,
When these systems of orthogonal
polynomials are moreover finite (so-called Leonard pairs) then our theory
is related to Terwilliger \cite{8}.
In section \ref{sec:08} we apply the general theory of section \ref{sec:06}
to the little $q$-Jacobi and the $q$-Hahn case.
Finally, in section \ref{sec:04} we give $p$-adic group interpretations of
the $q=0$ results obtained in section \ref{sec:01}.

In an earlier version \cite{32} of this paper we also discussed
big $q$-Jacobi polynomials. However, this family does not fit nicely into the
general theory of this paper, and interpretations on $p$-adic groups are
yet missing.
\paragraph{Acknowledgement}
We thank Erik Koelink, Michael Schlosser, Paul Terwilliger, Michael Voit and
the referee for helpful comments.
\section{Little $q$-Jacobi polynomials}
\label{sec:01}
%until 01.44
%
In this section we describe the main themes of this paper through
the little $q$-Jacobi polynomials. These polynomials and their
$p$-adic limit were the pioneering example which motivated this
paper. The cellular structure of the $p$-adic Hecke algebra of
Grassmannians of lines (see \S\ref{sec:04}) is essentially an $LU$
factorization of the matrix corresponding to the spherical functions
w.r.t.\ a geometric basis. As this basis is the limit of little
$q$-Jacobi polynomials, it was natural to seek for such
factorization for the matrix corresponding to the latter
polynomials.

Throughout we assume $0<q<1$. See standard formulas for
little $q$-Jacobi polynomials in \cite[\S3.12]{1}.
\subsection{Limit for $q\downarrow0$}
{\em Little $q$-Jacobi polynomials} are given by
\beq
p_n(x;a,b;q):=
\qhyp21{q^{-n},abq^{n+1}}{aq}{q,qx}\qquad(n\in\Znonneg).
\label{eq:01.01}
\eeq
For $b=0$ they
are known as {\em little $q$-Laguerre polynomials} or {\em Wall polynomials}
$p_n(x;a;q):=p_n(x;a,0;q)$, see \cite[\S3.20]{1}.
It will turn out that we have to rescale the parameters $a$ and $b$ in order
to be able to take the limit of these polynomials for $q\downarrow0$.
We define:
\begin{align}
p_n^{a,b;q}(x):=
p_n(q^x;q^{-1}a,q^{-1}b;q)=
\qhyp21{q^{-n},abq^{n-1}}a{q,q^{x+1}}\qquad&
\label{eq:01.02}
\\ \nonu
(n\in\Znonneg,\;x\in\Znonneg\union\{\iy\},\;0<a<1,\;b<1).&
\end{align}

By \cite[(3.12.2)]{1} the functions \eqref{eq:01.02}
satisfy the orthogonality relation
\beq
\sum_{x=0}^\iy p_m^{a,b;q}(x)\,p_n^{a,b;q}(x)\,w_x^{a,b;q}=
\frac{\de_{m,n}}{\om_n^{a,b;q}}\qquad(m,n\in\Znonneg),
\label{eq:01.03}
\eeq
where
\bea
w_x^{a,b;q}&:=&
\frac{(a;q)_\iy}{(ab;q)_\iy}\,\frac{(b;q)_x}{(q;q)_x}\,a^x,
\label{eq:01.04}
\\
\om_n^{a,b;q}&:=&\frac{1-abq^{2n-1}}{1-abq^{n-1}}\,
\frac{(a,ab;q)_n}{a^n\,(q,b;q)_n}\,.
\label{eq:01.05}
\eea
Note that the weights $w_x^{a,b;q}$ and the dual weights
$\om_n^{a,b;q}$ are positive under the constraints for $a$ and $b$ given
in \eqref{eq:01.02}.
Since the little $q$-Jacobi are orthogonal polynomials with respect
to an orthogonality measure of bounded support, they form a complete
orthogonal system in the $L^2$ space with respect to this measure,
so the functions \eqref{eq:01.02} also form a complete orthogonal
system in $\ell^2(\Znonneg; w^{a,b;q})$.

As was observed in \cite[(4.1), (4.2)]{6},
\cite[Remark 3.1]{19} and
\cite[(5.1), (5.3)]{5}, little $q$-Jacobi
polynomials can alternatively be expressed as a terminating
${}_3\phi_1$ by
application of the transformation formula
\cite[(III.8) or Exercise 1.15~(ii)]{2}
to the terminating ${}_2\phi_1$ in
\eqref{eq:01.02}. We obtain:
\beq
p_n^{a,b;q}(x)=q^{\half n(n-1)}\,(-a)^n\,\frac{(b;q)_n}{(a;q)_n}\,
\qhyp31{q^{-n},abq^{n-1},q^{-x}}b{q,\frac{q^{x+1}}a}.
\label{eq:01.20}
\eeq
As was also observed in the papers just quoted,
this is related to the fact that little $q$-Jacobi polynomials are
the duals of {\em $q^{-1}$-Al-Salam Chihara polynomials}
\bea
&&Q_n\left(\thalf(aq^{-x}+a^{-1}q^x);a,b\mid q^{-1}\right)
\nonu\\
&&\qquad\qquad=
(-1)^n b^n q^{-\half n(n-1)}\,((ab)^{-1};q)_n\,
\qhyp31{q^{-n},q^{-x},a^{-2}q^x}{(ab)^{-1}}{q,q^n ab^{-1}}.
\label{eq:01:42}
\eea
{}From \eqref{eq:01.20} and \eqref{eq:01:42} we get indeed the duality
\bea
&&(-a)^{-n}\,q^{-\half n(n-1)}\,\frac{(a;q)_n}{(b;q)_n}\,p_n^{a,b;q}(x)
=
\frac{(-1)^x(q^{-1}a^{-1}b)^{\half x} q^{\half x(x-1)}}{(b;q)_x}\,\nonu\\
&&\qquad\qquad
\times Q_x\left(\thalf(qab)^{-\half}(q^{-n}+abq^{n-1});
(qab)^{-\half},(qab^{-1})^\half;q^{-1}\right)\quad
(n,x\in\Znonneg)\quad
\label{eq:01.43}
\eea

Formula \eqref{eq:01.20} can be rewritten as:
\beq
p_n^{a,b;q}(x)=\sum_{k=0}^{\min(n,x)}
q^{\half(n-k)(n-k-1)}\,(-a)^{n-k}\,
\frac{(bq^k;q)_{n-k}\,(abq^{n-1};q)_k}{(a;q)_n\,(q;q)_k}\,
(q^n,q^x;q^{-1})_k\quad(x\in\Znonneg).
\label{eq:01.21}
\eeq
We obtain as an immediate corollary of \eqref{eq:01.21}:
\begin{theorem} \label{th:01.15}
The limit functions {\em (little 0-Jacobi functions)}
\beq
p_n^{a,b;0}(x):=\lim_{q\downarrow0}\,p_n^{a,b;q}(x)\quad
(x\in\Znonneg)
\eeq
exist. They are equal to
\bea
p_0^{a,b;0}(x)&=&\qquad\quad\; 1,
\label{eq:01.12}
\mLP
p_1^{a,b;0}(x)&=&
\left\{\begin{array}{cl}\dstyle -\,\frac{a(1-b)}{1-a}&
\mbox{if $x=0$,}\bLP
1&\mbox{if $x>0$,}
\end{array}\right.
\label{eq:01.13}
\mLP
p_n^{a,b;0}(x)&=&
\left\{\begin{array}{cl}
\quad 0&\;\mbox{\quad if $0\le x<n-1$,}\mLP
\dstyle -\,\frac a{1-a}&\;
\mbox{\quad if $x=n-1$,}\bLP
\quad 1&\;\mbox{\quad if $x>n-1$}
\end{array}\right\}
\quad(n\ge2).
\label{eq:01.14}
\eea
\end{theorem}

\begin{remark} \label{th:01.27}
Theorem \ref{th:01.15} was first stated by Haran
\cite[(7.3.37)]{3}, where the limit functions
\eqref{eq:01.12}--\eqref{eq:01.14} are given in \cite[(4.4.9)]{3}.
The limit result there follows from the expression for little
$q$-Jacobi polynomials on \cite[p.59, second formula from
below]{3}, which reads in our notation as:
\bea
p_n^{a,b;q}(x)&=&q^{nx}\,\frac{(b;q)_n}{(a^{-1}q^{1-n};q)_n}\,
\qhyp32{q^{-n},a^{-1}q^{1-n},q^{-x}}{b,0}{q,q}
\label{eq:01.28}
\\
&=&
\sum_{k=0}^{\min(n,x)}
q^{\half(n-k)(n+2x-3k-1)}\,
\frac{(-a)^{n-k}\,(bq^k;q)_{n-k}\,}{(a;q)_{n-k}\,(q;q)_k}\,
(q^n;q^{-1})_k\,(q^x;q^{-1})_k\,.\quad
\label{eq:01.29}
\eea
Formula \eqref{eq:01.28} follows from \eqref{eq:01.02} by the
transformation formula \cite[(1.5.6)]{2}.
Theorem \ref{th:01.15} can be obtained from \eqref{eq:01.29}
 by letting $q\downarrow0$.

Haran's \cite[Ch.~7]{3} notation is connected with ours by:
\beq \frac{\phi_{q,n}^{(\al)\be}(g^x)}{\phi_{q,n}^{(\al)\be}(0)}=
p_n^{q^\be,q^\al;q}(x), \qquad
\ze_{(q)}(s)=\left((q^s;q)_\iy\right)^{-1}. \label{eq:01.34} \eeq

Haran \cite[pp.\ 61--64]{3}
also considers the little $q$-Laguerre case ($\al\to\iy$ in
\eqref{eq:01.34}; $b=0$ everywhere in our notation)
and its $q=0$ limit.
\end{remark}
\begin{remark}
{}From \eqref{eq:01.21} we also obtain the following asymptotics of
$p_n^{a,b;q}(x)$ as $q\downarrow0$:
\beq
\lim_{q\downarrow0}q^{-\half (n-x)(n-x-1)}\,
p_n^{a,b;q}(x)=\frac{(-a)^{n-x}}{1-a}\qquad(1\le x\le n-1).
\label{eq:01.09}
\eeq
Alternatively, \eqref{eq:01.09} can be derived from the $q$-difference equation
\cite[(3.12.5)]{1} by induction with respect to $x$, starting at $x=1$.
Theorem~\ref{th:01.15} can also be proved by use of \eqref{eq:01.09}.
\end{remark}
{}From  \eqref{eq:01.04} and \eqref{eq:01.05} we get limits
\bea
w_x^{a,b;0}:=\lim_{q\downarrow0}\,w_x^{a,b;q}&=&
\left\{\begin{array}{cl}\dstyle \frac{1-a}{1-ab}&\;
\mbox{if $x=0$,}\bLP
\dstyle\frac{(1-a)(1-b)}{1-ab}\,a^x&\;\mbox{if $x>0$,}
\end{array}\right.
\label{eq:01.17}
\mLP
\om_n^{a,b;0}:=\lim_{q\downarrow0}\,\om_n^{a,b;q}&=&
\left\{\begin{array}{cl}
1&\quad\mbox{if $n=0$,}\mLP
\dstyle \frac{1-a}{a(1-b)}&\quad
\mbox{if $n=1$,}\bLP
\dstyle\frac{(1-a)(1-ab)}{a^n(1-b)}&\quad
\mbox{if $n>1$.}
\end{array}\right.
\label{eq:01.18}
\eea
Note that $w_x^{a,b;0}>0$ ($x\in\Znonneg$) and
$\om_n^{a,b;0}>0$ ($n\in\Znonneg$) if
$0<a<1$ and $b<1$.
The orthogonality relation \eqref{eq:01.03}
remains valid for $q=0$, as can be verified by use of
\eqref{eq:01.12}--\eqref{eq:01.14}
and \eqref{eq:01.17}--\eqref{eq:01.18}.
Formally, we can obtain the case $q=0$ of
\eqref{eq:01.03} by taking termwise limits.
\subsection{$LU$ factorization}
\label{sec:07}
Formula \eqref{eq:01.21}
has the big advantage over
\eqref{eq:01.29} that it
can be rewritten in matrix form as a product of a lower and an upper triangular
matrix:
\beq
P^q=L^q\,U^q,\quad{\rm i.e.,}\quad
P_{n,x}^q=\sum_{k=0}^\iy L_{n,k}^q\,U_{k,x}^q
=\sum_{k=0}^{\min(n,x)} L_{n,k}^q\,U_{k,x}^q\qquad(n,x\in\Znonneg),
\label{eq:01.23}
\eeq
where
\bea
P_{n,x}^q:=p_n^{a,b;q}(x),\;&&
L_{n,k}^q:=q^{\half(n-k)(n-k-1)}\,(-a)^{n-k}\,
\frac{(bq^k;q)_{n-k}\,(abq^{n-1};q)_k}{(a;q)_n\,(q;q)_k}\,
(q^n;q^{-1})_k\,,
\nonu\\
&&U_{k,x}^q:=(q^x;q^{-1})_k\,.
\label{eq:01.24}
\eea

As a limit case of \eqref{eq:01.23}, \eqref{eq:01.24}
formulas \eqref{eq:01.12}--\eqref{eq:01.14}
can similarly be rewritten in matrix form \eqref{eq:01.23} with $q=0$,
where
\bea
P_{n,x}^0&:=&p_n^{a,b;0}(x),
\nonu\\
U_{k,x}^0&:=&\left\{\begin{array}{cl}
1&\mbox{if $k\le x$,}\\
0&\mbox{if $k>x$,}
\end{array}\right.
\nonu\\
L_{n,k}^0&:=&\quad\,0\quad\mbox{if $\;k\ne n\;$ or $\;n-1$,}
\nonu\\
L_{n,n-1}^0&:=&
\left\{\begin{array}{cl}\dstyle -\,\frac {a(1-b)}{1-a}&\;
\mbox{if $n=1$,}\bLP
\dstyle -\,\frac a{1-a}&\;
\mbox{if $n>1$,}
\end{array}\right.
\nonu\mLP
L_{n,n}^0&:=&
\left\{\begin{array}{cl}
\; \quad1&\quad
\mbox{if $n=0$,}\mLP
\quad\dstyle\frac {1-ab}{1-a}&\quad
\mbox{if $n=1$,}\bLP
\quad\dstyle\frac 1{1-a}&\quad
\mbox{if $n>1$.}
\end{array}\right.
\label{eq:01.26}
\eea
Again, $L^0$ is a lower triangular and $U^0$ an upper triangular matrix.
The functions $c_k^q(q^x):=U_{k,x}^q/U_{k,k}^q=(q^x;q^{-1})_k/(q;q)_k$
and $c_k^0(x):=U_{k,x}^0$
(see \eqref{eq:01.24} and \eqref{eq:01.26}) can be considered as
forming a {\em cellular basis} in the terminology of
\cite[\S3.3]{4} and \cite{18}. The functions $q^x\to\const U_{k,x}^q$
can also be considered as the one-variable cases of
Okounkov's \cite{21}
{\em shifted
Macdonald polynomials}.
\begin{remark} \label{th:6.33}
For $0\le q<1$ we can consider $P^q$ as the matrix of a unitary operator from
the Hilbert space $\ell^2(\Znonneg;(w^q)^{-1})$ onto the Hilbert space
$\ell^2(\Znonneg;\om^q)$. Let $v_k:=a^{-k}$.
Probably, $U^q$ is the matrix of
a bounded linear operator from $\ell^2(\Znonneg;(w^q)^{-1})$ to
$\ell^2(\Znonneg;v)$, and $L^q$ is the matrix of
a bounded linear operator from $\ell^2(\Znonneg;v)$ to
$\ell^2(\Znonneg;\om^q)$.
Probably, the operators corresponding to $U^q$ and $L^q$ have bounded
inverses with matrices given by the explicit matrix inverses of $U^q$ and
$L^q$.
\end{remark}
\subsection{A product formula for little 0-Jacobi functions}
\begin{theorem}
\label{th:01.30}
The functions $p_n^{a,b;0}$ satisfy the product formula
\bea
p_n^{a,b;0}(x)\,p_n^{a,b;0}(y)=
\sum_{z=0}^\iy c_{x,y,z}^{a,b,0}\,p_n^{a,b;0}(z),
\label{eq:01.31}
\eea
where $c_{x,y,z}^{a,b,0}$
is given by
\beq
c_{x,y,z}^{a,b,0}=\left\{
\begin{array}{cl}
\de_{z,\min(x,y)}&\mbox{if $x\ne y$,}\mLP
\dstyle\frac{1-2a+ab}{1-a}&\mbox{if $x=y=z=0$,}\bLP
(1-b)a^z&\mbox{if $0=x=y<z$,}\mLP
0&\mbox{if $x=y>z$,}\mLP
\dstyle\frac{1-2a}{1-a}&\mbox{if $x=y=z>0$,}\bLP
a^{z-x}&\mbox{if $0<x=y<z$.}
\end{array}
\right.
\label{eq:01.33}
\eeq
In particular,
\beq
p_n^{a,b;0}(x)\,p_n^{a,b;0}(y)=p_n^{a,b;0}(\min(x,y))\quad
\mbox{if $x\ne y$.}
\eeq
Under the constraints $0<a<1$, $b<1$ we have $c_{x,y,z}^{a,b,0}\ge0$ for all
$x,y,z\in\Znonneg$ iff
\beq
0<a\le\thalf\quad\mbox{\rm and}\quad 2-a^{-1}\le b<1.
\label{eq:01.32}
\eeq
\end{theorem}
\noindent
{\bf Proof}\quad
Straightforward verification by
\eqref{eq:01.12}--\eqref{eq:01.14}.\qed
\begin{corollary}
The functions $p_n^{a,b;0}$ satisfy the product formula
\bea
p_n^{a,b;0}(x)\,p_n^{a,b;0}(y)=\frac{1-ab}{1-a}\,
\sum_{z=0}^\iy C_{x,y,z}^{a,b,0}\,w_z^{a,b,0}\,p_n^{a,b;0}(z),
\eea
where
\beq
C_{x,y,z}^{a,b,0}=\frac{1-a}{1-ab}\,
\sum_{n=0}^{\min(x,y,z)+1}
p_n^{a,b,0}(x)\,p_n^{a,b,0}(y)\,p_n^{a,b,0}(z)\,\om_n^{a,b,0}.
\eeq
is symmetric in $x,y,z$, and
for $x\le y\le z$ explicitly given by
\beq
C_{x,y,z}^{a,b,0}=\left\{
\begin{array}{cl}
0&\mbox{if $0\le x<y\le z$,}\mLP
1&\mbox{if $0=x=y<z$,}\mLP
(1-b)^{-1}a^{-x}&\mbox{if $0<x=y<z$,}\mLP
\dstyle\frac{1-2a+ab}{1-a}&\mbox{if $0=x=y=z$,}\bLP
\dstyle\frac{1-2a}{(1-a)(1-b)}\,a^{-x}&\mbox{if $0<x=y=z$.}
\end{array}
\right.
\eeq
\end{corollary}
\begin{remark} \label{th:01.44}
Dunkl and Ramirez \cite{23} obtained the little
$0$-Laguerre functions $p_n^{a,b;0}$ for $(a,b)=(p^{-1},0)$ as
spherical functions on the ring of $p$-adic integers, and they derived
the above product formula for those special parameter values from that
interpretation as spherical functions.
\end{remark}
\section{$LU$ factorizations: the general case}
\label{sec:06}
%until 06.36
%
\subsection{Lower times upper}
Let us put the results of \S\ref{sec:07}
in a more general framework.
Let $\FSN:=\{0,1,\ldots,N\}$ or $\Znonneg$
and let
$Y:=\{y_x\}_{x\in\FSN}$ be a
countable subset of $\RR$. Let $\{p_n\}_{n\in\FSN}$ be a complete system
of orthogonal polynomials on $Y$ with respect to positive weights $w_x$
on the points $y_x$:
\beq
\label{eq:06.01}
\sum_{x\in \FSN}p_n(y_x)\,p_m(y_x)\,w_x=(\om_n)^{-1}\,\de_{n,m}\quad
(n,m\in\FSN),
\eeq
where $\om_n>0$ for all $n\in\FSN$ and where, by completeness, we
also have the dual orthogonality relation
\beq
\label{eq:06.02}
\sum_{n\in\FSN}p_n(y_x)\,p_n(y_{x'})\,\om_n=(w_x)^{-1}\,\de_{x,x'}\quad
(x,x'\in \FSN).
\eeq
Note that completeness will certainly hold if $Y$ is finite or
bounded in  $\RR$.

We will now introduce some square matrices with row and column indices
running over $\FSN$.
Let $P$ be the matrix with entries $P_{n,x}:=p_n(y_x)$ ($n,x\in\FSN$).
Let also $W$ be the diagonal matrix with diagonal entries $w_x$
($x\in\FSN$) and let $\Om$ be the diagonal matrix with diagonal entries
$\om_n$ ($n\in\FSN$).
Then \eqref{eq:06.01} and \eqref{eq:06.02} can be written in matrix form as,
respectively,
\beq
PWP^t=\Om^{-1},\quad
P^t\Om P=W^{-1}.
\label{eq:06.16}
\eeq
Define polynomials $c_k$ of degree $k$ by
\beq
\label{eq:06.04}
c_k(y):=\prod_{j=0}^{k-1} (y_j-y)\quad(k\in\FSN,\;y\in\RR),
\eeq
and let $C$ be the matrix with entries
given by
\beq
\label{eq:06.03}
C_{k,x}:=c_k(y_x)=\prod_{j=0}^{k-1} (y_j-y_x)\quad(k,x\in\FSN).
\eeq
Then $C_{k,x}=0$ if $k>x$, so $C$ is an upper triangular matrix.
Then, for certain unique coefficients $B_{n,k}$ with $B_{n,n}\ne0$
we have:
\bea
p_n(y)&=&\sum_{k=0}^n B_{n,k} c_k(y)\quad(n\in\FSN,\;y\in\RR),
\label{eq:06.05}
\mLP
P_{n,x}&=&\sum_{k=0}^{\min(n,x)} B_{n,k} C_{k,x}\quad(n,x\in\FSN),
\label{eq:06.06}
\mLP
P&=&BC,
\label{eq:06.07}
\eea
where $B$ is the lower triangular matrix corresponding to the
coefficients $B_{n,k}$ ($n\ge k$).

Both $B$ and $C$ have two-sided inverses because they are triangular
matrices with nonzero diagonal entries. Furthermore, if $P=B'C'$
is another factorization of $P$ with $B'$ lower triangular and
$C'$ upper triangular, then $B'=BD$, $C'=D^{-1}C$ for some invertible
diagonal matrix $D$.

Things become even nicer if we know that there exist
orthogonal polynomials dual to $\{p_n\}$, i.e., if there exist
polynomials $r_x$ of degree $x$ ($x\in\FSN$) and a subset
$\{z_n\}_{n\in\FSN}$ of $\RR$ such that
\beq
\label{eq:06.08}
p_n(y_x)=r_x(z_n)\quad(n,x\in\FSN).
\eeq
Then the polynomials $r_x$ ($x\in\FSN$) will form a complete
system of orthogonal polynomials with respect to the weights
$\om_n$ on the points $z_n$ ($n\in\FSN$). Thus we can apply
the previous result to the $r_x$. Put
\beq
\label{eq:06.09}
B_{n,k}:=\prod_{i=0}^{k-1}(z_i-z_n)\quad(n,k\in\FSN).
\eeq
Then $B$ is a lower triangular matrix and for some upper triangular
matrix $C'$ we have
\beq
\label{eq:06.10}
p_n(y_x)=r_x(z_n)=\sum_{k=0}^{\min(n,x)} B_{n,k}C_{k,x}'\quad
(n,x\in\FSN).
\eeq
Then
$P=BC'$. Hence
\beq
P=BDC
\label{eq:06.18}
\eeq
for some diagonal matrix $D$
with nonzero diagonal entries~$\de_k$. Hence
\beq
\label{eq:06.11}
p_n(y_x)=\sum_{k=0}^{\min(n,x)}\de_k
\prod_{i=0}^{k-1}(z_i-z_n)\,
\prod_{j=0}^{k-1}(y_j-y_x).
\eeq
\subsection{Inverting the matrices $B$ and $C$ and computing $\de_k$}
The following theorem is the special case $a_j=1$, $b_j=0$ of
the Theorem in \cite[p.48]{20}, and it is also the case
$f(x):=(x_1-x)\ldots(x_{m-1}-x)$ ($m\le n$) of
\cite[p.54, Exercise 97]{22} (we thank Michael Schlosser for this reference),
but we will give here an independent proof.
\begin{theorem} \label{th:06.19}
For distinct complex numbers $y_n$ ($n\in\FSN$)
let $C=(C_{m,n})_{m,n\in\FSN}$ be an upper triangular matrix given by
\eqref{eq:06.03}.
Then
\beq
(C^{-1})_{k,n}=\prod_{\bisub{j=0}{j\ne k}}^n(y_j-y_k)^{-1}
\quad(0\le k\le n).
\label{eq:06.35}
\eeq
\end{theorem}
\noindent
{\bf Proof}\quad
Let $m\le n$. Put $V_{m,n}:=\prod_{m\le i<j\le n}(y_j-y_i)$.
We have to show that
\beq
\sum_{k=m}^n
\prod_{j=0}^{m-1}(y_j-y_k)\,
\prod_{\bisub{j=0}{j\ne k}}^n
(y_j-y_k)^{-1}=\de_{m,n}.
\label{eq:6.15}
\eeq
This is clearly true for $m=n$. For $m<n$
the \LHS\ of \eqref{eq:6.15} can be rewritten as
\begin{eqnarray*}
&&\sum_{k=m}^n
\prod_{\bisub{j=m}{j\ne k}}^n(y_j-y_k)^{-1}
\\
&&=\left(V_{m,n}\right)^{-1}\sum_{k=m}^n(-1)^{k-m}
\prod_{\bisub{m\le i<j\le n}{i,j\ne k}}(y_j-y_i)
\\
&&=\left(V_{m,n}\right)^{-1}
\sum_{k=m}^n(-1)^{k-n}
\left|\begin{matrix}
1&\ldots&1&1&\ldots&1\\
y_m&\ldots&y_{k-1}&y_{k+1}&\ldots&y_n\\
\vdots&&\vdots&\vdots&&\vdots\\
y_m^{n-m-1}&\ldots&y_{k-1}^{n-m-1}&y_{k+1}^{n-m-1}&\ldots&y_n^{n-m-1}
\end{matrix}\right|
\\
&&=\left(V_{m,n}\right)^{-1}
\left|\begin{matrix}
1&\ldots&1\\
1&\ldots&1\\
y_m&\ldots&y_n\\
\vdots&&\vdots\\
y_m^{n-m-1}&\ldots&y_n^{n-m-1}
\end{matrix}\right|=0.\qed
\end{eqnarray*}

It follows by taking transpose in \eqref{eq:06.35} that the inverse of the
lower triangular matrix $B$ given by \eqref{eq:06.09} is the lower
triangular matrix $B^{-1}$ with
\beq
(B^{-1})_{m,n}=\prod_{\bisub{i=0}{i\ne n}}^m(z_i-z_n)^{-1}\qquad
(0\le n\le m).
\label{eq:06.36}
\eeq

Put $y_\nu:=y_N$ if $\FSN=\{0,1,\ldots,N\}$ and put
$y_\nu:=\lim_{x\to\iy}y_x$ if $\FSN=\Znonneg$ and if
the (finite) limit $y_\nu$ exists and is not equal to any $y_x$ ($x\in\FSN$).
We will now derive an explicit expression for the coefficients $\de_k$
in \eqref{eq:06.11} involving $p_n(y_\nu)$ for all $n\in\FSN$.
This expression will have practical usage if a simple explicit expression
for $p_n(y_\nu)$ is known, as is the case in most examples.

It follows from \eqref{eq:06.11} that
\[
p_n(y_\nu)=\sum_{k=0}^n\de_k
\prod_{i=0}^{k-1}(z_i-z_n)\,
\prod_{j=0}^{k-1}(y_j-y_\nu)=
\sum_{k=0}^nB_{n,k}\,\de_k\prod_{j=0}^{k-1}(y_j-y_\nu).
\]
Hence, by matrix inversion,
\[
\sum_{n=0}^m\left(B^{-1}\right)_{m,n} p_n(y_\nu)=
\de_m\prod_{j=0}^{m-1}(y_j-y_\nu).
\]
Thus, by \eqref{eq:06.36} we obtain the following formula for $\de_m$:
\beq
\de_m=\frac1{\prod_{j=0}^{m-1}(y_j-y_\nu)}\,
\sum_{n=0}^m p_n(y_\nu)\prod_{\bisub{i=0}{i\ne n}}^m(z_i-z_n)^{-1}.
\label{eq:06.23}
\eeq
\subsection{Upper times lower}
We obtain from \eqref{eq:06.16} that
$P=\Om^{-1}(P^t)^{-1}W^{-1}$
and from
\eqref{eq:06.18} (only formally in the infinite dimensional case)
that $(P^t)^{-1}=(B^t)^{-1} D^{-1} (C^t)^{-1}$. Hence,
\beq
P=\Om^{-1}(B^t)^{-1}D^{-1}(C^{-1})^tW^{-1}.
\label{eq:06.20}
\eeq
{}From Theorem \ref{th:06.19},\eqref{eq:06.09} and
\eqref{eq:06.03} we see that
\bea
\left((B^t)^{-1}\right)_{n,k}&=&
\prod_{\bisub{i=0}{i\ne n}}^k(z_i-z_n)^{-1}\qquad
(0\le n\le k),
\label{eq:06.22}
\\
\left(C^{-1}\right)_{x,k}&=&
\prod_{\bisub{j=0}{j\ne x}}^k(y_j-y_x)^{-1}\qquad
(0\le x\le k).
\eea
When we substitute everything in \eqref{eq:06.20}
then we obtain
\beq
p_n(y_x)=\om_n^{-1}w_x^{-1}
\sum_{k\ge\max(n,x)}\de_k^{-1}
\prod_{\bisub{i=0}{i\ne n}}^k(z_i-z_n)^{-1}
\prod_{\bisub{j=0}{j\ne x}}^k(y_j-y_x)^{-1}.
\label{eq:06.21}
\eeq
\begin{remark} \label{th:6.34}
The derivation of \eqref{eq:06.20} is purely formal if $\FSN=\Znonneg$,
since we do not know in general if the matrices $B$ and $C$ correspond to
bounded linear operators and if these operators have bounded inverses.
See \cite{28} and \cite{27} and references given there for
some generalities about existence of $LU$-factorizations of bounded
linear operators as a product of a lower triangular and an upper triangular
matrix, both corresponding to bounded linear operators.
See \cite{29} for an example of a unitary operator on $\ell^2(\Znonneg)$
without $LU$-factorization.
\end{remark}
\section{$LU$ factorizations: examples}
\label{sec:08}
%until 08.04
%
\subsection{Little $q$-Jacobi}
{}From \eqref{eq:01.02} and \eqref{eq:01.43} we see that
\eqref{eq:06.11} will have meaning with the following substitutions:
\beq
y_x:=q^x,\quad
z_i:=q^{-i}+abq^{i-1}\quad
(x,i\in\Znonneg),\qquad
\frac{p_n(y_x)}{p_n(y_\nu)}:=p_n^{a,b;q}(x),
\label{eq:08.01}
\eeq
where
\beq
y_\nu=\lim_{x\to\iy}y_x=0,\qquad
p_n(y_\nu)=(-a)^{-n}\,q^{-\half n(n-1)}\,\frac{(a;q)_n}{(b;q)_n}.
\label{eq:08.02}
\eeq
Rewrite \eqref{eq:06.11} as
\beq
\frac{p_n(y_x)}{p_n(y_\nu)}=\sum_{k=0}^{\min(n,x)}
\left(p_n(y_\nu)\right)^{-1}
\left({\textstyle\prod_{i=0}^{k-1}(z_i-z_n)}\right)\,
\de_k\left({\textstyle\prod_{j=0}^{k-1}y_j}\right)\,
\left({\textstyle\prod_{j=0}^{k-1}\bigl(1-y_j^{-1}y_x\bigr)}\right)
\eeq
and compare with \eqref{eq:01.23}, \eqref{eq:01.24}.
We obtain that
\[
U_{k,x}^q=\textstyle\prod_{j=0}^{k-1}\bigl(1-y_j^{-1}y_x\bigr),\quad
L_{n,k}^q=\left(p_n(y_\nu)\right)^{-1}
\left({\textstyle\prod_{i=0}^{k-1}(z_i-z_n)}\right)\,
\de_k\left({\textstyle\prod_{j=0}^{k-1}y_j}\right).
\]
Hence, by \eqref{eq:01.24}, \eqref{eq:08.01} and \eqref{eq:08.02},
\beq
\de_m=\frac{(qa^{-1})^m}{(b;q)_m\,(q;q)_m}\,.
\label{eq:06.13}
\eeq

We can alternatively compute $\de_m$ from \eqref{eq:06.23}, which takes after
substitution of \eqref{eq:08.01} and \eqref{eq:08.02} the form
\beq
\de_m=\frac{q^m}{(ab,q;q)_m}
\qhyp64{q^{-1}ab, (qab)^\half,-(qab)^\half,a,0,q^{-m}}
{(q^{-1}ab)^\half,-(q^{-1}ab)^\half,b,abq^m}
{q,q^{-1}a^m}.
\label{eq:06.24}
\eeq
The ${}_6\phi_4$ can be evaluated as a confluent limit case of the
summation formula for a terminating very well-poised ${}_6\phi_5$ series
given in \cite[(2.4.2)]{2}. The resulting explicit formula for $\de_m$
by \eqref{eq:06.24} coincides with \eqref{eq:06.13}.

Note that \eqref{eq:06.11} with substitution of \eqref{eq:08.01} and
\eqref{eq:08.02} does not immediately allow to take limits for
$q\downarrow0$. For this we have to renormalize $p_n(y_x)$ by division by
$p_n(y_\nu)$ and we have to transfer some factors only depending on $k$
from $C_{k,x}$ to $B_{n,k}\de_k/p_n(y_\nu)$.
\mPP
Next we consider formula \eqref{eq:06.21}
(the upper times lower factorization)
with substitution of \eqref{eq:08.01} and
\eqref{eq:08.02}. Then $w_x$ and $\om_n$ in \eqref{eq:06.21} become
\[
w_x=w_x^{a,b;q},\qquad
\om_n:=(p_n(y_\nu))^2\,\om_n^{a,b;q}.
\]
With these substitutions, \eqref{eq:06.21} can be written in the form
\bea
p_n^{a,b;q}(x)&=&
\sum_{k=\max(n,x)}^\iy
q^{\half(k-x)(k-x-1)}\,(-a)^{k-x}
\frac{(abq^{n+k};q)_\iy\,(b;q)_k}{(a;q)_\iy\,(b;q)_x\,(q;q)_k}\,
(q^k;q^{-1})_n\,(q^k;q^{-1})_x\qquad
\label{eq:01.38}
\\
&=&\left\{\begin{array}{ll}
q^{\half(n-x)(n-x-1)}\,(-a)^{n-x}\,
\dstyle\frac{(abq^{2n};q)_\iy\,(b;q)_n\,(q^n;q^{-1})_x}
{(a;q)_\iy\,(b;q)_x}&\\
\qquad\times\dstyle\qhyp22{q^{n+1},bq^n}{q^{n-x+1},abq^{2n}}{q,q^{n-x}a}
&\mbox{if $n\ge x$,}
\bLP
\dstyle\frac{(abq^{n+x};q)_\iy\,(q^x;q^{-1})_n}{(a;q)_\iy}\,
\qhyp22{bq^x,q^{x+1}}{q^{x-n+1},abq^{n+x}}{q,a}&
\mbox{if $n\le x$.}
\end{array}\right.
\label{eq:01.39}
\eea

Above we started with \eqref{eq:06.21}. Instead we might have started with
\eqref{eq:06.20} and then obtain for $C$ and $B$
(\eqref{eq:06.03} and \eqref{eq:06.09} with substitution of
\eqref{eq:08.01}) the inverse matrices by
\cite[(4.2) and (4.11)]{7}.

Since \eqref{eq:06.21} was only derived in a formal way, we have not yet
proved now \eqref{eq:01.38} and \eqref{eq:01.39} in a rigorous way.
However,
the expressions \eqref{eq:01.39} can be alternatively
obtained from \eqref{eq:01.20}
by first applying \cite[(III.8) or Exercise 1.15~(ii)]{2} and
next inverting the series (see \cite[Exercise 1.4(ii)]{2}):
\bea
&&p_n^{a,b;q}(x)=
q^{\half n(n-1)}\,(-a)^n\,\frac{(bq^x;q)_n}{(a;q)_n}\,
\qhyp21{q^{-n},q^{-x}}{b^{-1}q^{1-n-x}}{q,q^{2-n}(ab)^{-1}}
\label{eq:01.40}
\\
&&\qquad=\left\{\begin{array}{ll}
q^{\half(n-x)(n-x-1)}\,(-a)^{n-x}\,
\dstyle\frac{(b;q)_n\,(q^n;q^{-1})_x}{(a;q)_n\,(b;q)_x}\,
\qhyp21{bq^n,q^{-x}}{q^{n-x+1}}{q,aq^n}&
\mbox{if $n\ge x$,}
\bLP
\dstyle\frac{(q^x;q^{-1})_n}{(a;q)_n}\,
\qhyp21{bq^x,q^{-n}}{q^{x-n+1}}{q,aq^n}&
\mbox{if $n\le x$.}
\end{array}\right.\qquad
\label{eq:01.41}
\eea
Then the expressions \eqref{eq:01.41} yield expressions \eqref{eq:01.39}
by means of \cite[(1.5.4)]{2}.
\subsection{$q$-Hahn}
See standard formulas for {\em $q$-Hahn polynomials} in
\cite[\S3.6]{1}.
They are given by
\beq
Q_n(x;a,b,N;q):=
\qhyp32{q^{-n},abq^{n+1},x}{aq,q^{-N}}{q,q}
\qquad(N\in\Zpos,\;n\in\{0,1,\ldots,N\}).
\label{eq:05.01}
\eeq
It will turn out that we have to renormalize the parameters $a$ and $b$
in order to be able to take the limit of these polynomials for $q\downarrow0$.
For the same reason, we have to consider these polynomials with argument
$q^{x-N}$, rather than the usual argument $q^{-x}$. We define:
\begin{align}
Q_n^{a,b,N;q}(x):=Q_n(q^{x-N};q^{-1}a,q^{-1}b,N;q)=
\qhyp32{q^{-n},abq^{n-1},q^{x-N}}{a,q^{-N}}{q,q}\qquad&
\label{eq:05.02}
\\ \nonu
\mbox{($N\in\Zpos$, $n,x\in\{0,1,\ldots,N\}$;
$0<a<1$, $b<1$ or $a,b>q^{1-N}$).}&
\end{align}
By \cite[(3.6.2)]{1} the functions \eqref{eq:05.02}
satisfy the orthogonality relation
\beq
\sum_{x=0}^N Q_m^{a,b,N;q}(x)\,Q_n^{a,b,N;q}(x)\,w_x^{a,b,N;q}=
\frac{\de_{m,n}}{\om_n^{a,b,N;q}}\qquad(m,n=0,1,\ldots,N),
\label{eq:05.08}
\eeq
where
\bea
w_x^{a,b,N;q}&:=&
\frac{(a;q)_N}{(ab;q)_N}\,
\frac{(q^N;q^{-1})_x}{(aq^{N-1};q^{-1})_x}
\frac{(b;q)_x}{(q;q)_x}\,a^x,
\label{eq:05.09}
\\
\om_n^{a,b,N;q}&:=&
\frac{(q^N;q^{-1})_n}{(abq^N;q)_n}\,
\frac{1-abq^{2n-1}}{1-abq^{n-1}}\,
\frac{(a,ab;q)_n}{a^n\,(q,b;q)_n}\,.
\label{eq:05.10}
\eea
Note that the weights $w_x^{a,b;q}$ and the dual weights
$\om_n^{a,b;q}$ are positive under the constraints for $a$ and $b$ given
in \eqref{eq:05.02}.

If we apply the transformation formula
\cite[(3.2.2)]{2}
to the terminating ${}_3\phi_2$ in \eqref{eq:05.02} then we obtain:
\bea
q^{-\half n(n-1)}\,(-a)^{-n}\,\frac{(a;q)_n}{(b;q)_n}\,
Q_n^{a,b,N;q}(x)&=&
\qhyp32{q^{-n},abq^{n-1},q^{-x}}{b,q^{-N}}{q,\frac{q^{x-N+1}}a}
\label{eq:05.19}
\\
&=&
\qhyp32{q^x,q^n,a^{-1}b^{-1}q^{-n+1}}{b^{-1},q^N}{q^{-1},q^{-1}}.
\label{eq:05.22}
\\
&=&
Q_n^{b^{-1},a^{-1};q^{-1}}(N-x)
\label{eq:05.23}
\\
&=&
Q_n(q^x;qb^{-1},qa^{-1},N; q^{-1})\quad
\label{eq:05.21}
\\
&=&
R_x(q^n+a^{-1}b^{-1}q^{-n+1};qb^{-1},qa^{-1},N\mid q^{-1}),\qquad
\label{eq:05.20}
\eea
which is both a polynomial of degree $n$ in $q^x$ for $n=0,1,\ldots,N$
and a polynomial of degree $x$ in $q^{-n}+abq^{n-1}$ for $x=0,1,\ldots,N$.
In \eqref{eq:05.20} we have used the notation
for {\em dual $q^{-1}$-Hahn polynomials} (see \cite[\S3.7]{1}
for dual $q$-Hahn polynomials; exchange there in \cite[(3.7.1)]{1} $n$ and $x$,
replace $q$ by $q^{-1}$, and
next replace $\ga$ by $qb^{-1}$ and $\de$ by $qa^{-1}$).

Formula \eqref{eq:05.19} can be rewritten as:
\beq
Q_n^{a,b,N;q}(x)=\sum_{k=0}^{\min(n,x)}
q^{\half(n-k)(n-k-1)}\,(-a)^{n-k}\,
\frac{(bq^k;q)_{n-k}\,(abq^{n-1};q)_k}{(a;q)_n\,(q;q)_k}\,
\frac{(q^n;q^{-1})_k\,(q^x;q^{-1})_k}{(q^N;q^{-1})_k}\,.
\label{eq:05.04}
\eeq
In notation \eqref{eq:05.02} and \eqref{eq:01.02} the limit formula
\cite[(4.6.1)]{1} reads
\beq
\lim_{N\to\iy} Q_n^{a,b,N;q}(x)=p_n^{a,b;q}(x).
\label{eq:05.05}
\eeq
Expressions \eqref{eq:05.02}, \eqref{eq:05.19}, \eqref{eq:05.04} for $q$-Hahn
polynomials and expressions \eqref{eq:05.09} and \eqref{eq:05.10} for their
weights and dual weights
are very similar to expressions
\eqref{eq:01.02}, \eqref{eq:01.20}, \eqref{eq:01.21},
\eqref{eq:01.04}, \eqref{eq:01.05},
respectively,
for little $q$-Jacobi polynomials, and the $q$-Hahn expressions immediately
turn into their little $q$-Jacobi counterparts under the limit transition
\eqref{eq:05.05}.
\begin{remark} \label{th:05.18}
Formula \eqref{eq:05.02} together with \eqref{eq:05.19} is the $q$-Hahn case
of a more general identity observed by Terwilliger \cite{8}
related to so-called parameter arrays (see Theorem 4.1 (ii), Lemma 4.2
and Example 5.4 in \cite{8}). The limit for $N\to\iy$ of this
$q$-Hahn case is the little $q$-Jacobi equality of
\eqref{eq:01.02} and \eqref{eq:01.20}. This observation may have relevance
for the open problem raised in \cite[Problem 11.1]{8}.
\end{remark}

{}From \eqref{eq:05.02} and \eqref{eq:05.19} we see that
\eqref{eq:06.11} will have meaning with the following substitutions:
\beq
y_x:=q^x,\quad
z_i:=q^{-i}+abq^{i-1}\quad
(x,i\in\{0,1,\ldots,N\},\qquad
\frac{p_n(y_x)}{p_n(y_N)}:=Q_n^{a,b,N;q}(x),
\label{eq:08.03}
\eeq
where
\beq
p_n(y_\nu)=(-a)^{-n}\,q^{-\half n(n-1)}\,\frac{(a;q)_n}{(b;q)_n}.
\label{eq:08.04}
\eeq
With these substitutions, formula \eqref{eq:06.11} coincides with
formula \eqref{eq:05.04} if we put
\beq
\de_m=\frac{(qa^{-1})^m}{(b;q)_m\,(q;q)_m\,(q^N;q^{-1})_m}\,.
\label{eq:06.26}
\eeq

We can alternatively compute $\de_m$ from \eqref{eq:06.23}.
This yields $\de_m$ as the \RHS\ of \eqref{eq:06.24}
with an additional factor $1/(q^N;q^{-1})_m$, from which we obtain
\eqref{eq:06.26}.
\mPP
Formula \eqref{eq:06.21} can also be specified in the $q$-Hahn case.
Make substitutions as above and, furthermore, put
\[
w_x:=w_x^{a,b,N;q},\qquad
\om_n:=\bigl((p_n(y_N)\bigr)^2 \om_n^{a,b,N;q}.
\]
Then \eqref{eq:06.21} yields
\bea
&&Q_n^{a,b,N;q}(x)=
\sum_{k=\max(n,x)}^N
q^{\half(k-x)(k-x-1)}\,(-a)^{k-x}\nonu
\\
&&\qquad\times
\frac{(q^N;q^{-1})_k\,(abq^{n+k};q)_{N-k}\,(aq^{N-1};q^{-1})_x\,
(b;q)_k}{(a;q)_N\,(b;q)_x\,(q;q)_k}\,
\frac{(q^k;q^{-1})_n\,(q^k;q^{-1})_x}{(q^N;q^{-1})_n\,(q^N;q^{-1})_x}\,.
\label{eq:06.27}
\eea
Formula \eqref{eq:06.27} can also be directly reduced to
\eqref{eq:05.02} by first reversing the direction of summation in
\eqref{eq:06.27}: substitute $k=N-l$. Then we obtain
\beq
Q_n^{a,b,N;q}(x)=q^{\half(N-x)(N-x-1)}\,(-a)^{N-x}\,
\frac{(bq^x;q)_{N-x}}{(a;q)_{N-x}}\,
\qhyp32{q^{n-N},q^{x-N},a^{-1}b^{-1}q^{-N-n+1}}{b^{-1}q^{-N+1},q^{-N}}{q,q}.
\label{eq:06.28}
\eeq
Formula \eqref{eq:06.28} follows from \eqref{eq:05.02}
by the transformation formula \cite[(3.2.2)]{2}.
Note that \eqref{eq:05.19} was obtained from \eqref{eq:05.02}
by a different application of this transformation formula.
Also note that the summation reversion changed the
upper times lower formula \eqref{eq:06.27} into
the lower times upper formula \eqref{eq:06.28}.

The ${}_3\phi_2$ in \eqref{eq:06.28} can be written both as a
{\em Hahn polynomial} and a
{\em dual $q$-Hahn polynomial}: we rewrite \eqref{eq:06.28}
for $n,x=0,1,\ldots,N$ as
\bea
&&q^{-\half(N-x)(N-x-1)}\,(-a)^{x-N}\,
\frac{(a;q)_{N-x}}{(bq^x;q)_{N-x}}\,
Q_n^{a,b,N;q}(x)
\nonu\\
&&\qquad=
Q_{N-n}(q^{x-N};b^{-1}q^{-N},a^{-1}q^{-N},N;q)
\label{eq:06.30}
\\
&&\qquad=
R_{N-x}(q^{n-N}+a^{-1}b^{-1}q^{-N-n+1};b^{-1}q^{-N},a^{-1}q^{-N}\mid q).
\qquad
\label{eq:06.29}
\eea
Thus the \LHS\ of the above identities is both a polynomial of degree $N-x$ in
$q^{n-N}+a^{-1}b^{-1}q^{-N-n+1}$ and a polynomial of degree $N-n$ in
$q^{x-N}$.
\begin{remark}
Analogous to our observation for \eqref{eq:05.19}
(see Remark \ref{th:05.18}), formula
\eqref{eq:06.28} can also be obtained as a consequence of
Theorem 4.1 (ii), Lemma 4.2
and Example 5.5 in \cite{8}. Also observe that finite orthogonal
polynomial systems whose duals are also orthogonal polynomial systems,
the so-called {\em Leonard pairs}, were extensively studied by
Terwilliger, see for instance \cite{24}, \cite{25}, \cite{8},
\cite{26}. Associated with a Leonard pair
is a {\em split decomposition}, which gives rise to a {\em parameter array}.
Formula (10) in \cite{8}, which depends on the parameters from that
array, is essentially the same as our formula \eqref{eq:06.11}.
\end{remark}
\subsection{$0$-Hahn functions}
We obtain as an immediate corollary of \eqref{eq:05.04}:
\begin{theorem} \label{th:05.06}
The limit functions {\em($0$-Hahn functions)}
\beq
Q_n^{a,b,N;0}(x):=\lim_{q\downarrow0}\,Q_n^{a,b,N;q}(x)\qquad
(n,x=0,1,\ldots,N)
\eeq
exist. They are equal to the little $0$-Jacobi functions
$p_n^{a,b;0}(x)$ (see \eqref{eq:01.12}--\eqref{eq:01.14})
restricted to $x=0,1,\ldots,N$:
\beq
Q_n^{a,b,N;0}(x)=p_n^{a,b;0}(x)\qquad(x=0,1,\ldots,N).\qquad
(n,x=0,1,\ldots,N)
\label{eq:05.07}
\eeq
\end{theorem}

{}From  \eqref{eq:05.09} and \eqref{eq:05.10} we get limits
\bea
w_x^{a,b,N;0}:=\lim_{q\downarrow0}\,w_x^{a,b,N;q}&=&
\left\{\begin{array}{cl}\dstyle \frac{1-a}{1-ab}&\;
\mbox{if $x=0$,}\bLP
\dstyle\frac{(1-a)(1-b)}{1-ab}\,a^x&\;\mbox{if $0<x<N$,}\bLP
\dstyle\frac{1-b}{1-ab}\,a^N&\;\mbox{if $x=N$,}
\end{array}\right.
\label{eq:05.11}
\mLP
\om_n^{a,b,N;0}:=\lim_{q\downarrow0}\,\om_n^{a,b,N;q}&=&
\left\{\begin{array}{cl}
1&\quad\mbox{if $n=0$,}\mLP
\dstyle \frac{1-a}{a(1-b)}&\quad
\mbox{if $n=1$,}\bLP
\dstyle\frac{(1-a)(1-ab)}{a^n(1-b)}&\quad
\mbox{if $2\le n\le N$.}
\end{array}\right.
\label{eq:05.12}
\eea
Note that $w_x^{a,b,N;0}>0$ ($0\le x\le N$) and
$\om_n^{a,b,N;0}>0$ ($0\le n\le N$) if
$0<a<1$ and $b<1$. They are almost the same as the weights
\eqref{eq:01.17} and duals weights \eqref{eq:01.18} for little
$0$-Jacobi polynomials.
The orthogonality relation \eqref{eq:05.08}
remains valid for $q=0$ by taking limits for $q\downarrow0$.
\begin{remark} \label{th:05.13}
Theorem \ref{th:05.06} was first stated
by Haran \cite[(7.3.37)]{3}, where the limit functions
\eqref{eq:05.07}
are given in \cite[(4.4.10)]{3}. The limit result there follows
(although this is not explicitly stated)
from the expression \cite[(7.3.20), (7.3.21)]{3} for $q$-Hahn polynomials,
which reads in our notation as:
\bea
&&Q_n^{a,b,N;q}(x)=q^{nx}\,
\frac{(q^{N-x-n+1};q)_n}{(q^{N-n+1};q)_n}
\frac{(b;q)_n}{(a^{-1}q^{1-n};q)_n}\,
\qhyp32{q^{-n},a^{-1}q^{1-n},q^{-x}}{b,q^{N-x-n+1}}{q,q}
\label{eq:05.14}
\\
&&=
\sum_{k=0}^{\min(n,x)}
q^{\half(n-k)(n+2x-3k-1)}\,
\frac{(-a)^{n-k}\,(q^{N-x};q^{-1})_{n-k}\,(bq^k;q)_{n-k}}
{(q^N;q^{-1})_n\,(a;q)_{n-k}\,(q;q)_k}\,
(q^n;q^{-1})_k\,(q^x;q^{-1})_k\,.\quad\qquad
\label{eq:05.15}
\eea
Formula \eqref{eq:05.14} follows from \eqref{eq:05.19} by the
transformation formula \cite[(3.2.5)]{2}.
Theorem \ref{th:05.06} can be obtained from \eqref{eq:05.15}
by letting $q\downarrow0$.

Haran's \cite[Ch.~7]{3} notation is connected with ours by:
\beq
\frac{\phi_{q(N),n}^{(\al)\be}(x,N-x)}{\phi_{q(N),n}^{(\al)\be}(N,0)}=
Q_n^{q^\be,q^\al,N;q}(x), \qquad
\ze_{(q)}(s)=\left((q^s;q)_\iy\right)^{-1}. \label{eq:05.17} \eeq
\end{remark}
\section{Interpretation as spherical functions over $p$-adic spaces}
\label{sec:04}
%until 04.01
%
We mention for completeness the overall picture concerning our main
object of study, the little $q$-Jacobi polynomials. Let $\FF$ be a
local field. $\FF$ can be Archimedean ($\RR$, $\CC$) or
non-Archimedean, that is, either
a finite extension of the field $\QQ_p$ of
$p$-adic numbers or
the Laurent series over a finite field
(see \cite[Chapter 4]{30} for
details). For $\FF$ non-Archimedean, let $\FSO$ stand for the ring
of integers. Let $K_{\FF}$ be the maximal compact subgroup of
$GL(d,\FF)$: the orthogonal group $O_d$, the unitary group $U_d$ or
$GL(d,\FSO)$. The natural representation of $K_{\FF}$ arising from
its action on the projective space is given by
\[
\rho^{\FF}:K_{\FF} \longrightarrow
U\big(L^2(\PP^{d-1}_{\FF})\big), \qquad
[\rho^{\FF}(g)f](x)=f(g^{-1}x).
\]
This representation admits a multiplicity free decomposition into
irreducible representations:
\[
L^2(\PP^{d-1}_{\FF})=\bigoplus_{n \in \Znonneg}
\FSU_{n}^{\FF}.
\]
The label $\FF$ on the various objects here emphasizes the
dependence on the field. However, the point here is that the
decomposition does not depend on the field. Moreover, the
irreducibles occurring in the decomposition for fixed $n$
correspond to each other when we go through the various fields (cf.\ \cite{3},
\cite{4}). This correspondence is realized by the observation
that, for all $\FF$, the little $q$-Jacobi polynomial of degree $n$
has limits
which are spherical functions in $\FSU_{n}^{\FF}$ for all
$\FF$'s. The orthogonality measure of these limit functions is
the projection of the Haar measure from $K_{\FF}$ to the space
$\PP^{d-1}_{\FF} \times_{K_{\FF}} \PP^{d-1}_{\FF}$, on which the
spherical functions live.

It also turns out (cf.~\cite{4}) that this scheme could be
generalized to representations arising from the action of these
groups on Grassmannians.
\subsection{Interpretation of little 0-Jacobi functions}
Let $\FF$ be a $p$-adic field, $\FSO$ the ring of integers, $\wp$
the maximal ideal in $\FSO$, $p^r$ the cardinality of the residue
field $\FSO/\wp$ ($p$ a prime number),
and $p\FSO=\wp^e$ ($e$ the ramification index, see
again \cite[Chapter 4]{30} for details). We look at the
representation of $GL(d,\FSO)$, the maximal compact subgroup of
$GL(d,\FF)$, defined by
\[
\rho \colon GL(d,\FSO) \longrightarrow
B\left(L^2(\PP^{d-1}(\FF))\right),\quad [\rho(g)f](x)=f(g^{-1}x),
\]
arising from the action of $GL(d,\FSO)$ on $\PP^{d-1}(\FF)
\cong\PP^{d-1}(\FSO)$. Let $P_m$ stand for the intersection of a
standard maximal parabolic subgroup of type $(m,d-m)$ in
$GL(d,\FF)$ with $GL(d,\FSO)$. In particular, $\PP^{d-1}(\FSO)
\simeq GL(d,\FSO)/P_1$. When we look at $P_m$-invariants in the
representation then we have
\[
L^2\left(\PP^{d-1}(\FSO)\right)^{P_m}=L^2\bigl(P_m \backslash
GL(d,\FSO) / P_1\bigr).
\]
The group $GL(d,\FSO)$ acts on $\PP^{d-1}(\FSO)$ and hence on its
quotients $\PP^{d-1}(\FSO /\wp^k)$. Denote the stabilizer of
$(1:0:\cdots:0) \in \PP^{d-1}(\FSO /\wp^k)$ in $GL(d,\FSO)$ by
$P_1^{(k)}$. The space $\PP^{d-1}(\FSO /\wp^k)$ is the boundary of
a ball of radius $k$ in the rooted tree with root valency
$\frac{p^{rd}-1}{p^r-1}$ (the cardinality of the projective space over
the residue field) and remaining vertices of degree $p^{rd}+1$. The
orbits of $P_m$ on $\PP^{d-1}(\FSO /\wp^k)$ consist of $k+1$
points and the orbits of the limit space are parameterized by
$\Znonneg\union\{\iy\}$:
\bea P_m \backslash GL(n,\FSO) /P_1 = \coprod_{k=0}^{\infty} P_m
\left(\begin{smallmatrix} 1 & 0 & 0 & \cdots & 0 \\ 0 & 1 & 0 &
\cdots & 0\\ & \vdots   & & \vdots & \\ \wp^k & 0 & \cdots & 0&
1\end{smallmatrix}\right) P_1 \simeq \Znonneg\union\{\iy\}.
\eea
The projection $\mu_p$ of the Haar measure to the orbit space is given
(see \cite{3}) by
\beq \mu_p(\{k\})=w_k^{a,b;0}\quad
(k\in\Znonneg,\;a=p^{-r(d-m)},\;b=p^{-rm}), \qquad \mu_p(\{\iy\})=0,
\label{eq:04.01}
\eeq
where the weights $w_k^{a,b;0}$, explicitly given by
\eqref{eq:01.17}, are limits for $q\downarrow0$ of the weights for
the orthogonality of the little $q$-Jacobi polynomials. Thus, one
obtains an interpretation of these weights on a $p$-adic space if
$a,b$ are as in \eqref{eq:04.01}. Moreover, for these values of the
parameters, the $p$-adic spherical functions (i.e., fixed vectors
in the representation $\rho$ under $P_m$), are  the functions
$p_n^{a,b,0}$ given by \eqref{eq:01.12}--\eqref{eq:01.14}. Hence,
they are limits of little $q$-Jacobi polynomials, as was shown in
\cite{3}.
\subsection{Interpretation of little 0-Laguerre functions}
For the interpretation of the little $q$-Laguerre polynomials at $q=0$,
we look at the action of $GL(m,\FSO)$ on $\FSO^m$. The orbits of
this action are characterized by the minimal valuation of the
entries of a vector. Hence, the orbits are: $\FSO^m \setminus
\wp\FSO^m$, $\wp\FSO^m \setminus \wp^2\FSO^m$, $\cdots$, $\{0\}$.
The
measure of the orbits is precisely the orthogonality measure for
the $0$-Laguerre functions $p_n^{p^{-rm},0;0}$:
\[
\mu(\wp^{j}\FSO^m \setminus \wp^{j+1}\FSO^m)=(1-p^{-rm})p^{-rjm} \qquad
(j\in \Znonneg).
\]
Moreover, by taking
$GL(m,\FSO)$-invariants in the decomposition to irreducibles of the
representation $L^2(\FSO^m)$, we find the $0$-Laguerre functions
$p_n^{p^{-rm},0;0}$ as fixed vectors. In the special case $m=r=1$ (group
of $p$-adic units acting on the ring of $p$-adic integers) this interpretation
was already obtained by Dunkl and Ramirez \cite{23}.

The two pictures (Jacobi and Laguerre) are related in the
following manner. One can restrict the action of $GL(d,\FSO)$ on
$\PP^{d-1}(\FSO)$, to an action of its subgroup $GL(m,\FSO)$
embedded in the top left corner on the subspace
$\{(x:0:\cdots:0:1) | ~x\in \FSO^m\} \subset \PP^{d-1}(\FSO)$.
This action is clearly the same as the action $GL(m,\FSO)$ on
$\FSO^m$. In terms of the parameters $a$ and $b$, this restriction
amounts to setting $b=0$, thus ignoring the irrelevant part of the
space.
\subsection{Product formula --- $p$-adic}
In this subsection we derive the $p$-adic product formula. We assume
that $a=p^{-r(d-m)}$ and $b=p^{-rm}$.
Then the nonnegativity conditions \eqref{eq:01.32} for the
product formula \eqref{eq:01.31} are valid.
Let
\bea \nu(i)=\sum_{j \ge i}\mu_p(j)=\left\{\begin{array}{cl}
1 & \mbox{if $i=0$,}\mLP \dstyle\frac{1-b}{1-ab}\,a^i&\mbox{if $i>0$.}
\end{array}\right. \eea
Our first step is to look at the spherical functions with a
different normalization, which makes them idempotents in the
convolution algebra $L^1\bigl(P_m \backslash GL(d,\FSO) /
P_1\bigr)$. We also rewrite them in terms of the measure, rather
than in terms of $a$ and $b$.
\bea
e_0(x)=\frac{p_0^{a,b;0}(\iy)}{\|p_0^{a,b;0}\|^2}\,p_0^{a,b;0}(x)&=
&\qquad\quad\;
1, \mLP
e_n(x)=\frac{p_n^{a,b;0}(\iy)}{\|p_n^{a,b;0}\|^2}\,p_n^{a,b;0}(x)&=&
\left\{\begin{array}{cl} \quad 0&\;\mbox{\quad if $0\le
x<n-1$,}\mLP \dstyle -\frac{1}{\nu(n-1)}&\; \mbox{\quad if
$x=n-1$,}\bLP \quad \dstyle
\frac{1}{\nu(n)}-\frac{1}{\nu(n-1)}&\;\mbox{\quad if $x>n-1$}
\end{array}\right\}
\quad(n\ge1).\qquad
\eea
Let $\{c_i=\mathbf{1}_{\{i,\ldots,\infty\}}\}$ and
$\{g_i=\mathbf{1}_{\{i\}}\}$. For $i \ge 0$ we have:
\begin{align}
c_i&=\sum_{j \ge i} g_i &g_i&=c_i-c_{i+1} &\quad\\
c_i&=\nu(i)\sum_{j \le i}e_j
&e_i&=\frac{1}{\nu(i)}c_i-\frac{1}{\nu(i-1)}c_{i-1} & (c_{-1}:=0)
\end{align}
The multiplication in the algebra, which is defined by declaring
that the $e_i$'s are idempotents, is given by:
\begin{align}
&e_i \star e_j=\delta_{ij}e_i,\\
 &c_i \star c_j=\nu(\max\{i,j\})c_{\min\{i,j\}},\\
 &g_i \star g_j =\left\{\begin{array}{cl}
\mu(\max\{i,j\})g_{\min\{i,j\}}\qquad \qquad \quad&\quad\mbox{if
$i \ne j$,}\mLP (\mu(i)-\nu(i+1))g_i+\mu(i)\sum_{j > i}g_j&\quad
\mbox{if $i=j$.}
\end{array}\right.
\end{align}
At this point we restrict to the $p$-adic case. In particular we
use the fact that the $e_i$'s and $g_i$'s are dual bases in the
sense that the former is an idempotent basis for the convolution
product, and the latter is an idempotent basis for the pointwise
product. The spherical transform intertwines these products and
bases. It follows that the multiplication table for the $g_i$'s
(after normalizing) with respect to the convolution product is the
pointwise multiplication for the idempotents, giving the desired
product formula.
If we normalize the $g_i$'s to be orthonormal, by setting
$\hat{g}_i=\frac{1}{\mu(i)}g_i$, we get:
\bea
 \hat{g}_i \star \hat{g}_j =\left\{\begin{array}{cl}
\hat{g}_{\min\{i,j\}}\qquad \qquad \quad&\quad\mbox{if $i \ne
j$,}\mLP \big(1-\frac{\nu(i+1)}{\mu(i)}\big)\hat{g}_i+\sum_{j >
i}\frac{\mu(j)}{\mu(i)}\hat{g}_j&\quad \mbox{if $i=j$.}
\end{array}\right.
\eea
Which agrees with the $c_{x,y,z}^{a,b,0}$ in \eqref{eq:01.33}.
%
%:References

%
\vspace{\bigskipamount}
\begin{footnotesize}
\begin{quote}
Tom H. Koornwinder\\
Korteweg-de Vries Institute, University of
Amsterdam,\\
Plantage Muidergracht 24, 1018 TV Amsterdam, The Netherlands
\sLP
{\tt thk@science.uva.nl}
\bLP
Uri Onn\\
Einstein Institute of Mathematics, Edmond Safra Campus, Givat Ram,\\
The Hebrew University of Jerusalem, Jerusalem 91904, Israel
\sLP
{\tt urion@math.huji.ac.il}
\end{quote}
\end{footnotesize}
\end{document}